\newif\ifsmfart
\IfFileExists{smfart.cls}
   {\documentclass[12pt,english]{smfart}
    \IfFileExists{smfenum.sty}{\usepackage{smfenum}}{}
    \usepackage{bull}
    \smfarttrue}
{\message{^^J*** It would be better to typeset this file with smfart.cls ***^^J^^J}

\documentclass[12pt]{amsart}}

\usepackage{amssymb}
\usepackage{epsfig}
\usepackage{graphicx}
\numberwithin{equation}{section}

\advance\headheight 2pt
\calclayout
\allowdisplaybreaks[3]

\theoremstyle{plain}
\newtheorem{prop}[subsection]{Proposition}

\newtheorem{coro}[subsection]{Corollary}

\newtheorem{lemm}[subsection]{Lemma}

\theoremstyle{definition}

\theoremstyle{remark}
\newtheorem{rem}[subsection]{Remark}

\newcommand{\Z}{\mathbb Z}

\newcommand{\Ker}{{\rm Ker}}

\def\ra{\rightarrow}

\def\Map{{\rm Map}}
\def\Diff{{\rm Diff}}

\makeatother
\makeatletter

\author{Fedor Bogomolov}
\address{Courant Institute of Mathematical Sciences, N.Y.U. \\
 251 Mercer str. \\
 New York, NY 10012, U.S.A.}
\email{bogomolo@cims.nyu.edu}

\author{Yuri Tschinkel}
\address{Department of Mathematics \\
         Princeton University\\ 
         Fine Hall, Washington Road\\
         Princeton, NJ 08544-1000,  U.S.A.}
\email{ytschink@math.princeton.edu}

\title[Symplectic fourfolds]{Simple 
examples of symplectic fourfolds 
with exotic properties}

\begin{document} 

\date{\today}



\begin{abstract}
We construct examples of 
simply connected nonalgebraic symplectic fourfolds
with a prescribed number of nonintersecting 
symplectic curves with positive self-intersections. 
\end{abstract}


\maketitle

\tableofcontents

\section{Introduction}
\label{sect:introduction}

Projective varieties of complex dimension
2 are basic examples of symplectic fourfolds.
Of course, the class of algebraic surfaces  
is much smaller. There are many examples of 
symplectic nonalgebraic varieties with various
distinguishing properties (see 
\cite{gompf2}, \cite{fin-stern1},\cite{fin-stern2}, \cite{abk}).  
In  this note we focus on embedded curves in 
the algebraic versus symplectic category.
It is known that an algebraic surface $X$
satisfies the following properties:
\begin{itemize}
\item the intersection form on the subspace of 
homology generated by complex curves in $X$ is hyperbolic;
\item the fundamental group of any smooth ample
curve in $X$ surjects onto the 
fundamental group of the surface.
\end{itemize} 
We construct simple examples of symplectic fourfolds
violating both of these properties.

\

\

{\bf Acknowledgments.} The first author is grateful to Vik. Kulikov
and S. Nemirovski for useful discussions. 
We would like to thank the organizers of the 
``Monodromy''-conference in Moscow, June 2001. 
Both authors were partially supported by the NSF.

\section{Construction}
\label{sect:const}

Let $C$ be an orientable compact
Riemann surface of genus $g=g(C)>0$. 
Consider a smooth fibration $X\ra C$ such that 
every fiber $X_c$ (with $c\in C$) is a smooth
orientable compact Riemann surface of genus $g(X_c)>2$.  
We assume that the monodromy along every loop in 
the base  is represented by an orientable 
automorphism and that $X\ra C$ admits a smooth section.

\begin{lemm}
\label{lemm:s}
Such a fibration is symplectic.
\end{lemm}

\begin{proof}
This is a well known fact, but we sketch a proof for
completeness. 
We build a fiberwise nondegenerate form as follows:
choose a standard basis of loops $a_i,b_i$ of the 
fundamental group of the base $C$. 
Each loop defines a monodromy diffeomorphism of $X$
(modulo isotopy), denoted by the same letter.    
The group $\Diff_g(C)$ of orientation preserving 
diffeomorphisms of $C$ admits a contraction onto the
group of volume-preserving diffeomorphisms. 
Therefore, we can find representatives for $a_i,b_i$ 
preserving the volume.
By invariance, this defines a vertical nondegenerate 
{\em closed} form 
$\omega^0$ on the preimage of a neighborhood of
the basic loops. 
The complement is a disc and the 
boundary is isomorphic in a standard
way (by the existence of a section) to ${\mathbb S}^1\times X$. 
The monodromy diffeomorphism on the fiber is 
isotopic to the identity. We can choose the isotopy to 
be volume-preserving. This defines a closed 2-form 
$\omega$ and a closed 2-form 
\begin{equation}
\label{eqn:omega}
\omega_X:=\omega_C +\lambda \omega
\end{equation}
(here $\omega_C$ is a 2-form on the base $C$ and 
$\lambda$ is an arbitrary nonzero constant). 
The form $\omega_X$ defines the symplectic 
structure on the fibration $X\ra C$. 
\end{proof}

\begin{rem}
The proof gives a recipe how to construct
symplectic fibrations (see \cite{endo}). Any relation in the mapping
class group $\Map_g$ of the form
\begin{equation}
\label{eqn:ll}
\prod_{i=1}^{g} [g_i,g_i'] = 1
\end{equation}
defines a symplectic 
fibration of the above type over
a Riemann surface $C$. 

\end{rem}

We will need standard extensions 
of the mapping class group $\Map_g$. 
Denote by $\Map_g(n)$ the mapping class group of
a Riemann surface of genus $g$ with $n$ distinct labeled points.
Let $\Map_g{\langle n\rangle}$ 
be the mapping class group preserving a small disc at each labeled
point. It is a standard central $\Z^n$-extension of $\Map_g(n)$:
the kernel is generated by Dehn twists in the
neighborhood of each labeled point. More precisely, a
vector 
$$
(\ell_1,...,\ell_n)\in 
\Ker(\Map_g{\langle n\rangle}\ra \Map_g(n))
$$ 
has the following geometric interpretation.
Lift the relation (\ref{eqn:ll}) into $\Map_g(n)$. 
The obtained family of Riemann surfaces
has $n$ nonintersecting sections $s_i$. We can compute
the squares of these sections by considering 
the product 
$\prod_{i=1}^{g(C)} [ \tilde{g_i},\tilde{g}_i']$
(where $\tilde{g}$ is an arbitrary lifting of $g$ into  
$\Map_g,{\langle n\rangle}$); it is an element of the center
of $\Map_g{\langle n\rangle}$ and thus a vector 
$(\ell_1,..., \ell_n)\in \Z^n$. We have $s_i^2=\ell_i$. 
Notice that for $g>2$ the groups  
$\Map_g, \Map_g(n)$ and $\Map_g{\langle n\rangle}$
are equal to their commutator subgroup.

\begin{rem}
The group $H_2(\Map_g(0),\Z)$ is equal to $\Z$ for
$g\ge 3$. There exists a linear lower bound for
the genus of smooth curves realizing a given class in 
this $H_2$ (see \cite{endo}, for example). 
The bound follows from the relation between this class and 
the signature of the corresponding symplectic fourfold.
Bounds  of such type appeared previously in 
the context of nilpotent groups in \cite{bg}. 
\end{rem}

\begin{lemm}
\label{lemm:sect}
For any vector $(\ell_1,...,\ell_n)$ there exists a curve
$C$, a smooth fibration $X\ra C$ into curves 
$X_c$ of genus $g(X_c)\ge 3$
and a set of smooth sections $s_1,...,s_n$ of 
this fibration such that $s_i^2=\ell_i$. 
\end{lemm}

\begin{proof}
Every element in $\Map_g{\langle n\rangle}$
is representable as a product of commutators.
It suffices to represent the 
central element $(\ell_1,...,\ell_n)$. 
\end{proof}

\begin{coro}
\label{coro:ind}
There exists a symplectic fourfold $X$ and 
a smooth symplectic curve $D\subset X$
with $D^2>0$ such that the image of 
$\pi_1(D)$ in $\pi_1(X)$ has infinite index.
In particular, $X$ has no topological Lefschetz pencils
containing $D$ (or its multiples $rD$).     
\end{coro}

\begin{proof}
Assume that $\ell_i>0$ for all $i=1,...,n$. 
Take $X$ as in Lemma~\ref{lemm:sect}.
Changing the symplectic form $\omega_X$ in \ref{eqn:omega} 
(by making $\lambda$ sufficiently large) 
we can insure that all sections $s_i$ are symplectic.
Take $D$ to be one of these sections. 
\end{proof}

\begin{rem}
This corollary corrects the argument in Section 4 of 
\cite{bk}. 
\end{rem}

\begin{coro}
\label{coro:hyp}
There exists a symplectic fourfold $X$ such that
the intersection form on symplectic curves $D\subset X$ 
is not hyperbolic.
\end{coro}

\begin{proof}
For any $n>1$ and any vector
$(\ell_1,..., \ell_n)$ with positive $\ell_i$
we choose $X$ as in Lemma~\ref{lemm:sect}.  The
restriction of the intersection form to 
the sections $s_i$ is positive, contradicting
hyperbolicity.  
\end{proof}

\begin{rem}
Surfaces of such type 
can also be obtained from complex Kodaira fibrations
by reversing the orientation. 
Then the smooth complex curves which have 
negative normal bundle are turned into curves with
positive self-intersection. 
\end{rem}

It is well understood that 
symplectic geometry is, in a sense, more flexible 
or closer to differential geometry and topology than 
to algebraic geometry if we allow large fundamental groups.
Now we show how to modify the above examples to obtain
simply connected symplectic varieties with the same 
interesting properties.

\begin{prop}
\label{prop:ex}
For any $n>0$ there exists a 
simply connected symplectic fourfold
containing $n$ smooth symplectic nonintersecting 
curves with positive self-intersection.  
\end{prop}

\begin{proof}
Choose a vector $(\ell_1,...,\ell_n)\in \Z^n$ such that 
all $\ell >1$. Choose $X\ra C$ as in Lemma~\ref{lemm:sect}. 
By construction, every fiber of 
$X\ra C$ is a symplectic subvariety. 
Blow up (symplectically) 
the intersection points of a fiber $X_0$ with 
the sections $s_{i}$.  
The obtained symplectic variety $\hat{X}\ra X$ has 
a collection of nonintersecting 
symplectic subvarieties: proper transforms $\hat{s}_i$ 
of the sections and $\hat{X}_0$ of the fiber $X_0$. 
We have $\hat{s}_i^2=\ell_i-1>0$ and $\hat{X}_0^2=-n$.
Choose 
an algebraic simply connected surface $V_0$ containing 
a  curve $Z_0$ (a symplectic surface) with 
self-intersection $Z_0^2=n$. 
Choose a rational (algebraic) surface $V_1$ containing 
a smooth algebraic curve $Z_1$ of genus $g(C)$ with 
self-intersection $Z_1^2=-\hat{s}_1^2$. 
We glue $V_0$ and $V_1$ to 
$\hat{X}$ along $Z_0$ and $\hat{X}_0$, resp.     
$Z_1$ and $\hat{s}_1$. We denote the obtained
fourfold by $\tilde{X}^{\rm sing}$.  
By results of Gromov~\cite{gromov} and Gompf~\cite{gompf1}, 
we can smooth symplectically $\tilde{X}^{\rm sing}$
without changing the  symplectic form outside of
a small neighborhood of $\hat{X}_0$ and $\hat{s}_1$. 
The resulting smooth symplectic variety $\tilde{X}$  
still contains $n-1$ nonintersecting  
symplectic curves $\tilde{s}_i$ with 
positive self-intersection.

We claim that $\tilde{X}$ is simply connected. 
Clearly, the singular variety $\tilde{X}^{\rm sing}$
is simply connected 
(the rational surface $V_1$ kills the
$\pi_1(\hat{s}_1)$,
$V_0$ kills $\pi_1(\hat{X}_0)$ and 
$\pi_1(\tilde{X}^{\rm sing})$
is generated by these subgroups).  
Finally, the smoothing doesn't 
change the fundamental group. 
\end{proof}

\begin{rem}
The surfaces in Proposition~\ref{prop:ex}
and Lemmas~\ref{coro:ind} and \ref{coro:hyp}
were constructed in 
response to a question of Vik. Kulikov. 
He pointed out that small symplectic quasi-complex deformations
of (the graphs) of algebraic surfaces with normal intersections
as in \cite{gompf1} or, more generally, in \cite{bk},  
fail to produce quasi-complex embedded symplectic curves
of the above type. 

\

Our approach is similar to \cite{endo}, though the precise result
appears to be new. 
\end{rem}



\begin{thebibliography}{99}


\bibitem{abk}
J. Amor\'os, F. Bogomolov, L. Katzarkov, T.  Pantev, 
{\em Symplectic Lefschetz fibrations with arbitrary fundamental groups},
J. Differential Geom. {\bf 54} (2000), no. 3, 489--545.


\bibitem{bg}
J. Barge ,  E. Ghys, 
{\em Surfaces et cohomologie borne}, 
Invent. Math. {\bf 92} (1988), 509--526.




\bibitem{bk}
F. Bogomolov, L. Katzarkov,
{\em Symplectic four-manifolds and projective surfaces},
Symplectic, contact and low-dimensional
topology (Athens, GA, 1996), {\em Topology Appl.}
{\bf  88} (1998), no. 1-2, 79--109. 

\bibitem{endo} 
H. Endo, M. Korkmaz, D. Kotschick, B. Ozbagci, A. Stipsicz,
{\em 
Commutators, Lefschetz fibrations and the signatures of surface bundles},
{\tt alg-geom 0103176}, (2001).


\bibitem{fin-stern1}
R. Fintushel, R. Stern, 
{\em Symplectic surfaces in a fixed homology class},
J. Differential Geom. {\bf 52} (1999), no. 2, 203--222. 


\bibitem{fin-stern2}
R. Fintushel, R. Stern,
{\em The canonical class of a symplectic 4-manifold},
Turkish J. Math. {\bf 25} (2001), no. 1, 137--145. 





\bibitem{gompf1}
R. E. Gompf, 
{\em A new construction of symplectic manifolds},
Ann. of Math. (2) {\bf 142} (1995), no. 3, 527--595. 




\bibitem{gompf2}
R. E. Gompf,
{\em The topology of symplectic manifolds},
Turkish J. Math. {\bf 25} (2001), no. 1, 43--59. 



\bibitem{gromov}
M. Gromov, 
{\em Partial differential relations},
Springer-Verlag, Berlin,  (1986).




\end{thebibliography}
\end{document}